\magnification=\magstep1
\documentstyle{amsppt}

\parskip=\baselineskip
\baselineskip=1.1\baselineskip
\parindent=0pt
\loadmsbm
\UseAMSsymbols
\raggedbottom

\def\doct{\delta_{oct}}
\def\pt{\hbox{\it pt}}
\def\Vol{\hbox{vol}}
\def\sol{\operatorname{sol}}
\def\dih{\operatorname{dih}}
\def\vor{\operatorname{vor}}
\def\score{\sigma}

\def\rad{\operatorname{rad}}
\def\span{\operatorname{span}}
\def\qed{{\it q.e.d.\/}}
\def\minus{{\scriptscriptstyle{-}}}

\def\qed{{\hbox{}\nobreak\hfill\vrule height8pt width6pt 
     depth 0pt}\medskip}

\def\diag|#1|#2|{\vbox to #1in {\vskip.3in\centerline{\tt Diagram #2}\vss} }
\def\v{\hskip -3.5pt }
\def\gram|#1|#2|#3|{
        {
        \smallskip
        \hbox to \hsize
        {\hfill
        \vrule \vbox{ \hrule \vskip 6pt \centerline{\it Diagram #2}
         \vskip #1in %
             \includegraphics{#3}\hrule }
        \v\vrule\hfill
        }
\smallskip}}

\centerline{\bf Sphere Packings II}
\bigskip
\centerline{Thomas C. Hales}
\bigskip
\footnote""{\line{\it\hfill published in Discrete and Computational Geometry, 18:135-149, 1997}}

{\bf Abstract:}  An earlier paper describes a 
program to prove the Kepler conjecture on sphere packings.
This paper carries out the second step of that program.
A sphere packing leads to a decomposition of ${\Bbb R}^3$
into polyhedra.  The polyhedra are divided into two
classes.  The first class of polyhedra, called quasi-regular
tetrahedra, have density at most that of a regular tetrahedron.
The polyhedra in the remaining class
have density at most that of a regular
octahedron  (about $0.7209$).  
 
\bigskip

\centerline{\bf Section 1. Introduction}

This paper is a continuation of the first part of this
series \cite{4}.  The terminology and notation of this
paper are consistent with this earlier paper, and we refer
to results from that paper by prefixing the relevant
section numbers with `I'.

We review some definitions from \cite{4}.
Begin with
a packing of nonoverlapping spheres of radius 1 in Euclidean three-space.
The {\it density} of a packing is defined in \cite{1}.  It is defined
as a limit of the ratio of the
volume of the unit balls in a large region of space
to the volume of the large region.
The density of the packing may be improved by adding spheres until
there is no further room to do so.  The resulting packing is said to
be {\it saturated\/}.  

Every saturated packing gives rise to a decomposition of space into
simplices called the {\it Delaunay decomposition} \cite{8}. 
The vertices
of each Delaunay simplex are centers of spheres of the packing.
By the definition of the decomposition, none of the centers
of the spheres of the packing lie in the interior of the circumscribing
sphere of any Delaunay simplex.
We refer to the centers of the packing as {\it vertices\/}.
Vertices that come within $2.51$ of each other 
are called {\it close neighbors}.

The Delaunay decomposition is dual to the well-known Voronoi
decomposition.
If the vertices of the Delaunay simplices are in nondegenerate
position, two vertices are joined by an edge exactly when
the two corresponding Voronoi cells share a face, three vertices
form a face exactly when the three Voronoi cells share an edge, and
four vertices form a simplex exactly when the four corresponding
Voronoi cells share a vertex.  In other words, two vertices are
joined by an edge if they lie on a sphere that does not contain
any other of the vertices, and so forth (again assuming the vertices
to be in nondegenerate position).  

We say that the convex hull of four vertices
is a
{\it quasi-regular tetrahedron\/} (or simply a {\it tetrahedron})
if all four vertices are close neighbors of one another.
If the largest circumradius of the faces of a Delaunay simplex
is at most $\sqrt{2}$, we say that the simplex is {\it small}.
Suppose that we have a configuration of
six vertices in bijection with the vertices
of an octahedron with the property that two vertices are close
neighbors if and only if the corresponding vertices of the
octahedron are adjacent.
Suppose further that there is a unique diagonal of
length at most $2\sqrt{2}$.
In this case we call the convex hull of the six vertices 
a {\it quasi-regular octahedron\/} 
(or simply an {\it octahedron}).
A {\it Delaunay
star\/} is defined as the collection of all quasi-regular tetrahedra,
octahedra,
and Delaunay simplices that share a common vertex $v$.  

We assume that every simplex $S$ in this paper comes with a fixed
order on its edges, $1,\ldots,6$.
The order on the edges is to be arranged so that
the first, second, and third edges meet at a vertex.  We may
also assume that the edges numbered $i$ and $i+3$ are opposite
edges for $i=1,2,3$.  We define $S(y_1,\ldots,y_6)$ to be
the (ordered) simplex whose $i$th edge has length $y_i$.
If $S$ is a Delaunay simplex in a fixed Delaunay star, then
it has a distinguished vertex, the vertex common to all simplices
in the star.  In this situation, we assume that the edges
are numbered so that the first, second, and third edges meet
at the distinguished vertex.

A function, known as the {\it compression\/} $\Gamma(S)$, is defined on
the space of all Delaunay simplices.  Set
$\delta_{oct} = (-3\pi+12\arccos(1/\sqrt{3}))/\sqrt{8}\approx 0.720903$.
Let $S$ be a Delaunay simplex.  Let $B$ be the union of
four unit balls placed at each of the vertices of $S$.
Define the compression as
$$\Gamma(S) = -\delta_{oct} \Vol(S) + \Vol(S\cap B).$$
We extend the definition of compression to Delaunay stars $D^*$
by setting $\Gamma(D^*) = \sum\Gamma(S)$, with the sum running
over all the Delaunay simplices in the star.
We define a {\it point} (abbreviated $\pt$) to be
$\Gamma(S(2,2,2,2,2,2))\approx 0.0553736$.
The compression is often expressed as a multiple of $\pt$.

There are several other functions of a Delaunay simplex that will be
used.  
The {\it dihedral angle\/} $\dih(S)$ is defined
to be the dihedral angle of the simplex $S$ along the first
edge (with respect to the fixed order on the edges of $S$).
The {\it solid angle\/} (measured in steradians) at the vertex
joining the first, second, and third edges is denoted
$\sol(S)$.
Let $\rad(S)$ be the circumradius of the simplex $S$.
More generally, let $\rad(F)$ denote the circumradius of the
face of a simplex.
Let $\eta(a,b,c)$ denote the circumradius of a triangle with
edges $a$, $b$, $c$.
Explicit formulas for all these functions appear in I.8.

\bigskip
Fix
a Delaunay star $D^*$ about a vertex $v_0$, which we
take to be the origin, and we consider the unit sphere at $v_0$.
Let $v_1$ and  $v_2$ be vertices of $D^*$ such that
$v_0$, $v_1$, and $v_2$ are all close neighbors of one another.
We take the radial projections $p_i$
of $v_i$ to the unit sphere with center at the origin
and connect the points $p_1$ and $p_2$ by a geodesic
arc on the sphere.
We mark all such arcs on the unit sphere.
The closures of the
connected components of the complement of these arcs are
regions on the unit sphere,
called the {\it standard regions}.  We may remove the arcs
that do not bound one of the regions.
The resulting system of edges and regions will be
referred to as the {\it standard decomposition\/} of the unit sphere.

Let $C$ be the cone with vertex $v_0$ over one of the
standard regions.
The collection of the Delaunay simplices,
quasi-regular tetrahedra, and quasi-regular octahedra
of $D^*$ in $C$ (together with the distinguished vertex $v_0$)
will be called a {\it standard cluster}.  Each
Delaunay simplex in $D^*$ belongs to a unique standard cluster.

A real number, called the {\it score},
 will be attached to each cluster.  
Each star
receives a score by summing the scores for the clusters in the
star.   

The steps of the Kepler conjecture, as outlined in Part I, are

{

\def\ha{ \hangindent=20pt \hangafter=1\relax }
1. A proof that even if all standard
regions are triangular, the total score
is less than $8\,\pt$

\ha
2.  A proof that the standard clusters 
with more than
three sides 
score at most $0\,\pt$

\ha
3. A proof that if all of the 
standard regions are triangles or quadrilaterals,
then the total score is less than $8\,\pt$ (excluding the
case of pentagonal prisms)

\ha
4.  A proof that if some 
standard region has
more than four sides, then the
star scores less than $8\,\pt$

\ha
5.  A proof that pentagonal prisms score less than $8\,\pt$
 
}

The proof of the first step is complete. The other steps are
briefly discussed in Part I.
This paper establishes step 2.  Partial results have been
obtained for step 3 \cite{5}.
C.A. Rogers has shown that the density of a regular tetrahedron
is a bound on the density of packings in ${\Bbb R}^3$ \cite{8}.  The main
result of this paper may be interpreted as saying that the
density ($\delta_{oct}\approx 0.7209$) of a regular octahedron
is a bound on the density of the complement in ${\Bbb R}^3$
of the quasi-regular tetrahedra
in the packing.  

The score of a Delaunay star is
obtained by mixing Delaunay stars with the dual Voronoi cells.
Delaunay stars $D^*$
and the associated function $\Gamma$ 
behave much better than estimates of density by Voronoi cells,
provided each Delaunay simplex in the Delaunay star has a
small circumradius.  Unfortunately, $\Gamma(S)$
gives an increasingly poor bound on the density as the circumradius
of the Delaunay simplex $S$ increases.  
When the
circumradius of $S$ is greater than about $1.8$, it becomes
extremely difficult to  prove anything about sphere packings with the
function $\Gamma(S)$.  The score is
introduced to regularize the irregular behavior of $\Gamma(S)$.

Voronoi cells also present enormous difficulties.
The dodecahedron shows that a single Voronoi cell cannot lead
to a bound on the density of packings better than about
$0.75$.  This led
L. Fejes T\'oth to propose an approach to the Kepler conjecture
in which two layers of
Voronoi cells are considered: one central
Voronoi cell and a number of surrounding ones.   
Wu-Yi Hsiang has
made some progress in this direction, but there
remain many technical difficulties \cite{3}, \cite{7}.

The method of scoring in this paper seeks to combine the best
aspects of both approaches.  When the circumradius of a simplex
is small, we proceed as in Part I.  However, when the circumradius
of a simplex is large, we switch to Voronoi cells.  
Remarkably, these two approaches may be coherently
combined to give a meaningful score to Delaunay stars and,
by extension, a bound on the density of a packing.
The calculations of this paper suggest that this
hybrid approach to packings retains the best features of both
methods with no (foreseeable) negative consequences.

\bigskip
\centerline{\bf Section 2.  Some polyhedra}
\bigskip

Sometimes the tip of a Voronoi cell protrudes beyond the
face of a corresponding Delaunay simplex  (see Diagram 2.1.a).
This section describes a construction that amounts to slicing
off the protruding
tip of a Voronoi cell and reapportioning it among
the neighboring cells (see Diagram 2.1.b).

\gram|1.5|2.1|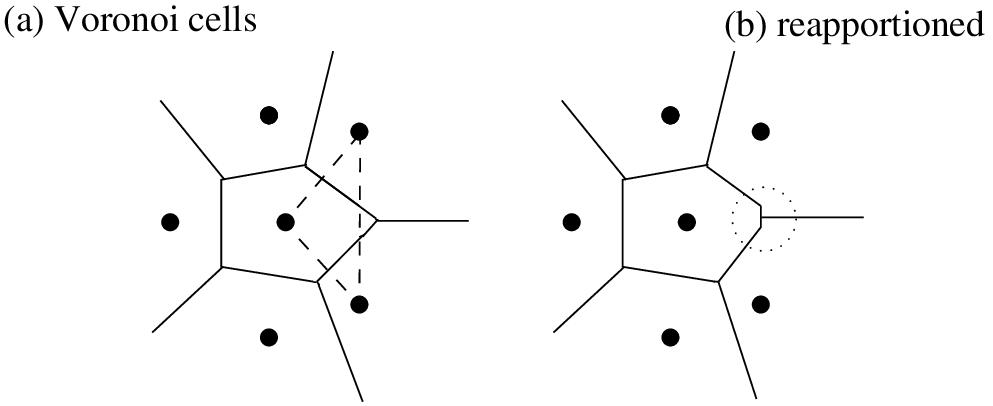|

Let $D^*$ be a Delaunay star with center $v_0=0$.  Let $V$ be the
Voronoi cell around $v_0$, obtained by duality from $D^*$.  
As a matter of convenience, we may assume that each point
in ${\Bbb R}^3$ belongs to a unique Voronoi cell by making
an arbitrary  choice for each point on the boundary of a cell.
If $R$ is a
standard cluster (possibly a single quasi-regular tetrahedron) in $D^*$, 
let $C(R)$ denote its cone over $v_0$:
$$C(R) = \{ t x : t\ge 0,\ x\in R\}.$$

In general, $V\cap C(R)$ depends
on more vertices than just those in the cluster $R$.
It is convenient to consider the slightly larger polyhedron
$V^0_R$ defined by just the vertices of $D^*$
that are in $R$.  That is, let $V^0_R$ be
the intersection of $C(R)$ with the half-spaces
$\{x : x\cdot v_i\le v_i\cdot v_i/2,\ \forall i\ne0\}$, 
where $\{v_i\}_i$ are the vertices
(other than $v_0$)
of the simplices and quasi-regular solids in the cluster $R$.
The faces of $V^0_R$ at $v_0$ are contained in the
triangular faces bounding the standard region of $R$.  The
other faces of $V^0_R$ are contained in planes
through the faces of the Voronoi cell $V$.  We refer
to these as Voronoi faces.
If $R$ is not a quasi-regular tetrahedron, set $V_R = V_R^0$.
If $R$ is a quasi-regular tetrahedron, we take the slightly smaller
polyhedron $V_R$ obtained by intersecting $V_R^0$
with the half-space (containing $v_0$)
bounded by the hyperplane through the
face of $R$ opposite the origin $v_0$.
(This may cut a tip from the Voronoi cell.)
By construction, $V_R$ depends only
on the simplices in $R$. 
The polyhedron $V_R$ is based at the center of some Delaunay star,
giving it a distinguished vertex $v$.  We write $V_R = V_R(v)$
when we wish to make this dependence explicit.  

By construction, $V_R^0 \supset V\cap C(R)$.
It is often true that $V_R = V\cap C(R)$.  Let us study the
conditions under which this can fail.  
We say that a vertex $w$ {\it clips\/} a standard cluster $R$
(based at $v_1$) if $w\ne v_1$ and some point of $V_R^0(v_1)$
belongs to the Voronoi cell at $w$.  Part I makes a 
thorough investigation of the geometry when a vertex
$w$ clips a quasi-regular tetrahedron.  (The vertex $w$ must
belong to a second quasi-regular tetrahedron that shares a 
face with $S$.  The shared face must have circumradius
greater than $\sqrt{2}$, and so forth.)

{\bf Lemma 2.2.}  Let $R$, based at a vertex $v_0$,
 be a standard cluster other than a quasi-regular tetrahedron.
Suppose it is clipped by a vertex $w$.  Then there is a face
$(v_0,v_1,v_2)$ of $R$ such that $(w,v_0,v_1,v_2)$ is a
quasi-regular tetrahedron.  Furthermore, $(v_0,v_1,v_2)$ is
the unique face of the quasi-regular tetrahedron of circumradius
at least $\sqrt{2}$.

\bigskip

{\bf Proof:}
Consider a point $p$ in $V_R \setminus V$.  Then there
exists a vertex $w\not\in C(R)$ of 
$D^*$ such that $p\cdot w > w\cdot w/2$.
The line segment from $p$ to $w$ intersects the cone $C(F)$ of
some triangular face $F$ that bounds the standard region of $R$
and has $v_0$ as a vertex.  Let $v_1$ and $v_2$ be the other vertices
of $F$.  By the construction of the faces bounding a standard region, the
edges of $F$ have lengths between $2$ and $2.51$.

Consider the region $X$ containing $p$ and bounded by the planes
$H_1 = \span(v_1,w)$, $H_2 = \span(v_2,w)$, $H_3 = \span(v_1,v_2)$,
$H_4 = \{x: x\cdot v_1 = v_1\cdot v_1/2\}$,  
and $H_5 = \{x: x\cdot v_2
= v_2\cdot v_2/2\}$.  The planes $H_4$ and $H_5$ contain the
faces of the Voronoi cell at $v_0$ defined by the vertices $v_1$
and $v_2$.  The plane $H_3$ contains the face $F$. The planes
$H_1$ and $H_2$ bound the region containing points, such as $p$,
that can be connected to $w$ by a segment that passes through $C(F)$.

Let $P = \{x: x\cdot w> w\cdot w/2\}$.  
The choice of $w$ implies
that $X\cap P$ is nonempty.  
We leave it as an exercise to check that $X\cap P$ is
bounded.  If the
intersection of a bounded polyhedron with a half-space is nonempty, then
some vertex of the polyhedron lies in the half-space.  So some
vertex of $X$ lies in $P$.

We claim that the vertex of $X$ lying in $P$ cannot lie on $H_1$.
To see this, pick coordinates $(x_1,x_2)$ on the plane $H_1$ 
with origin $v_0=0$ so that 
$v_1 = (0,z)$ (with $z>0$)
 and $X\cap H_1\subset X':=\{(x_1,x_2) : x_1\ge 0, \ x_2\le z/2\}$.
See Diagram 2.3.  If $X'$ meets $P$,
then the point $v_1/2$
lies in $P$.  This is impossible, because every point between
$v_0$ and $v_1$ lies in the Voronoi
cell at $v_0$ or $v_1$, and not in the Voronoi cell of $w$.  
(Recall that $|v_1-v_0|<2.51<2\sqrt{2}$.)

\bigskip
\gram|1.5|2.3|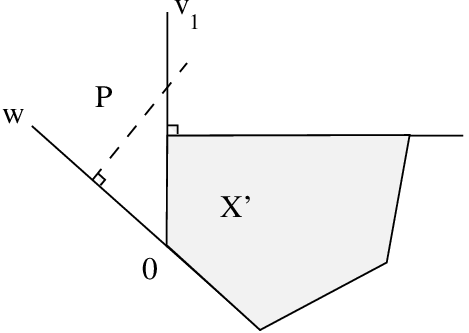|
\smallskip

Similarly, the vertex of $X$ in $P$ cannot lie on $H_2$.  Thus, the
vertex must be the unique vertex of $X$ that is not on $H_1$ or
$H_2$, namely, the point of intersection of $H_3$, $H_4$, and $H_5$.
This point is the circumcenter $c$ of the face $F$.  We conclude
that the polyhedron $X_0:= X\cap P$ contains $c$.
Since $c\in X_0$, the hypotheses of Lemma I.3.4 are met for $T=F$,
and
the vertices $v_0$, $v_1$, $v_2$, and $w$ 
are the vertices of a quasi-regular tetrahedron $S$. 
By I.3.4,
$|w'-v_i|< 2.3$, for $i=0,1,2$. The circumradius of the face
$F$ is between $\sqrt{2}$ and $2.51/\sqrt{3}\approx 1.449$.\qed

In the same context, if $w$ and $w'$ both clip $R$, then the
regions they cut from $V_R(v_0)$ are disjoint.  For
otherwise, a common point would belong to both
$V_S^0(w)$ and $V_S^0(w')$ where $S$ and $S'$ are the
two quasi-regular tetrahedra constructed by the Lemma.
Part I.3 shows that $S$ and $S'$ share their unique face
of circumradius greater than $\sqrt{2}$.  This is impossible,
because the lemma states that this face is shared with $R$.

Although the polyhedron $X_0$ 
belongs to the Voronoi cell
at $w$, it is included in the polyhedron $V_R$.  Similarly,
by repeating the construction at $v_1$ and $v_2$, we find that
there are small regions $X_1$, $X_2$ 
(with vertex $c$) in polyhedra $V_{R_1}$ and
$V_{R_2}$ at
$v_1$ and $v_2$ respectively 
that belong to the Voronoi cell at $w$.  

Call the union $X_0\cup X_1\cup X_2$ the {\it tip} protruding
from the quasi-regular tetrahedron $S$.  Associated
with a quasi-regular tetrahedron is at most one such tip. (The tip must
protrude from the face of $S$  with circumradius 
greater than $\sqrt{2}$.)
By construction, the tip is the set of points
$$\{ x: | x-w| \le |x-v_i|,\ \hbox { for } i=0,1,2;\ \ 
	\det(x,v_1,v_2)\det(w,v_1,v_2) \le 0\}.$$
This is $V_S^0(w)\setminus V_S(w)$.

The tip is a subset of the Voronoi cell at $w$.
Section I.3 explains the conditions under which this can fail to
hold.  There must be another vertex $u\ne w$ with the property
that $|u-v_i|< 2.3$, for $i=0,1,2$.  Then $u$, $v_0$, $v_1$, and $v_2$
are the vertices of a second quasi-regular tetrahedron $S'$ with face $F$, 
and this is contrary to our assumption that $R$ is not a quasi-regular
tetrahedron.
\bigskip

{\bf Corollary 2.4.}  The polyhedra $V_R$ cover ${\Bbb R}^3$ evenly as
we range over all the standard clusters of all the Delaunay stars
of the packing.

\bigskip
{\bf Proof:}
The preceding analysis shows that the polyhedra $V_R$
are obtained from the Voronoi cells by taking each protruding
tip, breaking it into three
pieces $X_0$, $X_1$, $X_2$, and attaching the piece $X_i$ to
the Voronoi cell at $v_i$.
The Voronoi cells cover ${\Bbb R}^3$ evenly.  As a result
of this analysis,
we see that
the polyhedra $V_R$ cover ${\Bbb R^3}$
evenly.\qed

To give one example of the size of the tip, we consider
the extreme case of the
tetrahedron $S=S(2,2,2,2.51,2.51,2.51)$.
Diagram 2.5 shows a correctly scaled drawing of
a tip protruding from the largest face
of $S$.
\bigskip

\gram|2.0|2.5|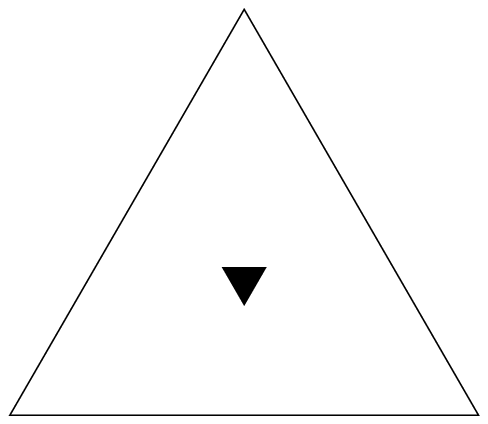|

\bigskip
\centerline{\bf Section 3.  The score attached to a Delaunay star}
\bigskip

This section gives some rules
for computing the score.  They were developed as a result of
computer experimentation suggesting when it is
advantageous to use Voronoi cells over Delaunay simplices.
This section actually gives an entire family of scoring systems.
This extra bit of flexibility will be useful as we encounter new
examples in the remaining steps of the program.
We expect the score to satisfy Conjecture I.2.2, which asserts that
the score of a Delaunay star is at most $8\,\pt$,
for all the scoring systems satisfying properties 1--4 below.
The Kepler conjecture is true if Conjecture I.2.2 holds for any one
such scoring system.
We have found through experimentation that
small or seemingly innocent changes in the score can lead to enormous
changes in the complexity of the optimization problem.  

{\bf 3.1.\ } This paper
proves the second step of the program for all of the scoring
systems presented below.
   Write $\sigma(S)$ for the score of $S$.

1.  Suppose that the standard cluster $R$ is a single
	quasi-regular tetrahedron: $R=S$.
	When the circumcenter of $S$ is contained in $S$,
	$$-4\delta_{oct} \Vol(V_R\cap C(S)) + 4 \sol(S)/3$$
	is an analytic function of the lengths of the
	edges.  This expression has an analytic continuation,
	denoted $\vor(S,V_R)$,
	to simplices $S$ that do not necessarily
	contain their circumcenter.

	If $\rad(S)>1.41$, then define the
	score to be $\vor(S,V_R)$.
	If $\rad(S)\le 1.41$, then define the
	score to be the compression $\Gamma(S)$.
	(This rule agrees with the definition 
	of $\vor(S)$ given in Section I.2.)

2.  Let $S$ be a small simplex that is not a quasi-regular tetrahedron.  
The score of $S$ will be either
$\vor(S)$ or $\Gamma(S)$ depending on criteria to be determined
by future research.\footnote"*"{More generally, we might add a small
constant c to the score of $S$ at one of its vertices and subtract
the same constant from another vertex.}
These criteria may depend on whether $S$ belongs to a quasi-regular
octahedron, but not on the position of any 
vertices of the packing outside $S$.
 It is essential for the scoring at all four vertices 
to have the same type (Voronoi or compression).
The only constraint imposed by the second step of the Kepler conjecture
will be $\sigma(S)\le 0$, if $S$ is small.  This leads to the
following mild restrictions on the use of Voronoi scoring.

If one of the first three edges is the long edge (say the first),
	compression scoring is to be used if the second, third, and
	fourth edges
	have length at most $2.06$, and the fifth and sixth edges have
	length at most $2.08$.

If one of the last three edges (say the fourth) is
	the long edge, compression scoring is to be used if
	(a), (b), (c), and (d) hold.

(a) The first edge has length at most at most $2.06$.

(b) The second and third edges have length at most $2.08$.

(c) The fifth and sixth edges have length at most $2.2$.

(d) The fourth edge has length at most $2.58$, or the
	fifth and sixth edges have lengths at most $2.12$.

3.  Suppose that $R$ is any standard cluster 
	other than a quasi-regular
	tetrahedron.  The cluster 
	is a union of Delaunay simplices
	$S_1,\ldots,S_r$.  Index the simplices so that 
	$S_1,\ldots,S_p$, for some $p\le r$ are the small simplices
	in the cluster.
	We define the score of the cluster $R$ to be
	$$\sum_{1\le i\le p} \sigma(S_i) +
		\sum_{p<i\le r} \vor(S_i,V_R),$$
where $\vor(S,V_R)= 4(-\doct\Vol(V_R\cap C(S))+\sol(S)/3)$.

4.  If $D^*$ is a Delaunay star, then its total score 
	$\sigma(D^*)$ is a sum of
	the scores of the standard clusters 
	of $D^*$.

\bigskip

Consider
the quasi-regular tetrahedron $S$
of Section 2 with vertices
$v_0$, $v_1$, $v_2$, and $w$ that has a protruding tip
$X_0\cup X_1\cup X_2$.  Let $\sol_v(S)$ denote the
solid angle of $S$ at the vertex $v$. The analytic continuation
$\vor(S,V_R)$ has the following geometric interpretation.
$$\vor(S,V_R(v)) = -4\delta_{oct}(\Vol(S,V_R(v)) + A(v)) + 
4\sol_v(S)/3,$$
with the correction term
$A(v_i) =  -\Vol(X_i)$, for  $i=0,1,2$,
and 
$$A(w)= \sum_{i=1}^3 \Vol(X_i).$$

The only pieces that are compression scored
are small simplices, everything else is Voronoi scored.  The small simplices
that are compression scored will be called simplices of {\it compression
type}.  The Voronoi-scored small simplices will be called simplices
of {\it Voronoi type}.
We define the {\it restricted cell\/} of a cluster $R$ to
be the complement in $V_R$ of the small simplices
in the packing.

\bigskip

{\bf Lemma 3.2.} (1) The score
of a cluster depends only on the cluster,
and not on the way it sits in a Delaunay
star or in the Delaunay decomposition of
space.

(2)  
Let $\Lambda$
denote the vertices of a saturated packing.  Let $\Lambda_N$
denote the vertices inside the ball of radius $N$. (Fix
any center for the ball.)  Let $D^*(v)$ denote the Delaunay star at
$v\in \Lambda$.  Then the score satisfies (in Landau's
notation)
$$\sum_{\Lambda_N} \score(D^*(v)) = \sum_{\Lambda_N} \Gamma(D^*(v)) + O(N
^2).$$

\bigskip

{\bf Proof:}
Statement (1) holds by construction.  

(2) The score reapportions the
compression of a given star among surrounding stars.
The second part of the lemma follows from the claim that 
everything is accounted
for, if we ignore the boundary effects caused by the truncation $N$.
Space is partitioned into regions each counted $-4\delta_{oct}$ times by
the compression of some star.  Each point in a sphere of the packing
is counted four times by the compression of some star.  
To verify 3.2.2, we must
check that the same holds of the score.

We switch from Voronoi to compression scoring on certain
small simplices.
The faces $F$ of a small simplex $S$ satisfy $\rad(F)\le\sqrt{2}$, so
no point on a face $F$ of $S$ can be closer to another vertex in the
packing than it is to the closest vertex of $F$.   This
has two implications.  First, 
the only polyhedra $V_R$ meeting a small 
simplex $S$ are the four based at the vertices of
$S$.  Second, Let $R$ be a standard cluster.
Let $S$ be a small simplex in $R$.  Then $V_R \cap C(S) = V_R \cap S$.
(In other words, tips cannot protrude from a small simplex.)
This means that 
the restricted cells and small simplices 
cover space evenly.  This decomposition is compatible
with the standard decomposition of a Delaunay star.

Consider the rules defining the score.  In counting the 
part of the volume
of a sphere contained in a simplex $S$, we see that it appears
four times with weight $1$ for a total weight of $4$,
	when $S$ is a small simplex of compression type.
It appears
once with weight $4$ for a total weight of $4$,
	when $S$ is of Voronoi type.

The result is now clear.\qed

{\bf Remark 3.3.}  It is useful to summarize the proof from
a slightly different point of view.  If $S$ is a 
quasi-regular tetrahedron or a small Delaunay simplex, then
the sum of its four scores, for each of its four vertices,
is $4\Gamma(S)$.  This follows directly from the definitions
(and the proof of Lemma 3.2) if the circumcenter of $S$ is contained in
$S$ (which is always the case for small simplices),
and it follows by analytic continuation in general.
Any other point in space belongs to a unique Voronoi cell centered
at some vertex $v$.  If the point is not in a tip protruding
from a quasi-regular tetrahedron, it is counted in the score at 
$v$.  If, however, the point belongs to a protruding tip, it
is counted in the score at exactly one of the three vertices,
other than $v$, of the quasi-regular tetrahedron.  In
this way, every point in ${\Bbb R}^3$ is accounted for.

{\bf Remark 3.4.}  The choice of the parameter $\mu=1.41$ in Rule 1
is somewhat arbitrary.  The choice is based on the comparison of
the functions
$$f_1(x) = \vor(S(2,2,2,2.51,2.51,x),V_S)\ \ \hbox{and }\ \ 
  f_2(x) = \Gamma(S(2,2,2,2.51,2.51,x)).$$
The difference $f_1(x)-f_2(x)$ has a zero for some $x\in [2.2603,2.2604]$.
This gives a crude estimate of when it is advantageous to
switch from $\Gamma(S)$ to $\vor(S,V_S)$.  The constant $1.41$
is a little more than $\rad(S(2,2,2,2.51,2.51,2.2604))\approx 1.405656$.

{\bf Proposition 3.5.}  The Delaunay stars in the face-centered
cubic and hexagonal-close packings score $8\,\pt$.

{\bf Proof:}  The eight regular tetrahedra each score $1\,\pt$, and 
 each regular octahedron scores $0\,\pt$, because
it has density $\delta_{oct}$,
for a total of $8\,\pt$.  \qed

We will see in Proposition 4.6 that the regular octahedron can be
broken into smaller pieces that score $0\,\pt$.

\bigskip
\centerline{\bf Section 4. The Main Theorem}
\bigskip

{\bf Theorem 4.1}  (a) The score of any small
quasi-regular tetrahedron is at most $1\,\pt$.
(b) The score of any other standard cluster is at most $0\,\pt$.

{\bf Proof.}  Statement (a) is a special case of Calculation I.9.1.
A quasi-regular
tetrahedron of Voronoi type scores less than $0\,\pt$ by Lemma I.9.17.
In the remainder of the proof, we actually prove a much
stronger statement.  We explicitly decompose each cluster
(other than a quasi-regular tetrahedron)
into a number
of pieces and show that the density of each piece is at most $\delta_{oct}$.
Since $\vor(S,V_R)$ and $\Gamma(S)$ are zero precisely when the
corresponding densities are $\delta_{oct}$ (or when
the volumes are zero), the theorem will follow.
The relevant pieces will be congruent to one of the following
{\it types}:

1.  A small simplex that is not a quasi-regular tetrahedron

2.  A set $\{t x: 0 \le t \le 1,\  x\in P_2\}\subset {\Bbb R}^3$,
where $P_2$ is a measurable set and every point of $P_2$ has
distance at least
$1.18$ from the origin (Diagram 4.2.a)

3.  A set $\{t x: 0\le t\le 1,\ x\in P_3\}\subset {\Bbb R}^3$,
where $P_3$ is a wedge of a disk of the form
$$P_3 = \{(x_1,x_2,x_3): x_3 = z_0,\ x_1^2+x_2^2\le 2,\ 0\le x_2\le
\alpha x_1\},$$
for some $\alpha>0$ and some $1\le z_0\le 1.18$
	(Diagram 4.2.b)

4.  A Rogers simplex $R(a,b,\sqrt{2})$  
	where $1\le a\le 1.18$ and $4/3\le b^2 \le 2$
(see Section I.8.6 and Diagram 4.2.c)

\bigskip
\gram|2.0|4.2|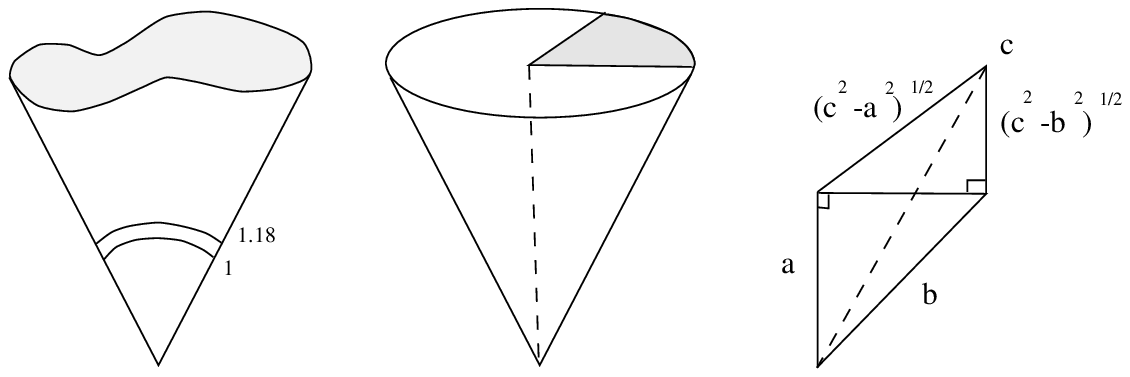|
\smallskip

In the first type, 
a unit ball is placed at each vertex of the simplex $S$,
and the density is the ratio of the volume of the 
part of the balls in $S$
to the volume of $S$.  In the second, third, and
fourth types, a unit ball is placed
at the origin, and the density is the ratio of the volume of the 
part of the ball
in the region to the volume of the region.

We decompose all of ${\Bbb R^3}$ into these four types and
quasi-regular tetrahedra. 
Set all the quasi-regular tetrahedra aside. 
Classify all
the small simplices, including those
contained in a quasi-regular octahedron, 
as regions of the first type. There remain the restricted cells.
Now fix a Delaunay star $D^*$, with center at the origin, and consider
the restricted cell of one of the clusters in the star.
We may assume that the restricted cell does not lie in a 
quasi-regular tetrahedron.
Break the restricted cell up further by taking its intersection
with the cones over each of its Voronoi faces $F$.  
Let $X$ be one such intersection.
If the face $F$ has distance
more than $1.18$ from the origin, classify $X$ 
as a region of the second
type.
Now assume the face $F$ has distance $h$ at most $1.18$ from the
center. 
Because $h<\sqrt{2}$, the point in the plane of $F$ closest
to the origin lies on the face $F$.
The set of points $P_2$ on the face $F$
at distance 
greater than $\sqrt{2}$ from the origin gives rise to
a region of the second type.  
To study what remains, we may truncate $F$
by intersecting it with a ball of radius $\sqrt{2}$.
Let $F'\subset F$ be the truncated face.

By Voronoi-Delaunay duality,
the face $F'$ lies in the bisecting plane between $0$ and some
vertex $v$ of the Delaunay star.
Consider the collection of triangles
formed by $0$, $v$, and another vertex of the Delaunay star $D^*$,
with the property that either
the triangle has circumradius 
at most $\sqrt{2}$ or all three
edges of the triangle have lengths between $2$ and $2.51$.
Consider the half-planes (bounded by the line through $0$ and $v$)
containing the various triangles in this collection.  This
fan of half-planes
partitions the face $F'$ into a collection of wedge-shaped pieces.
Consider one of them, $F''$.  
We claim that it has the form of Diagram 4.3.

\bigskip
\gram|1.8|4.3|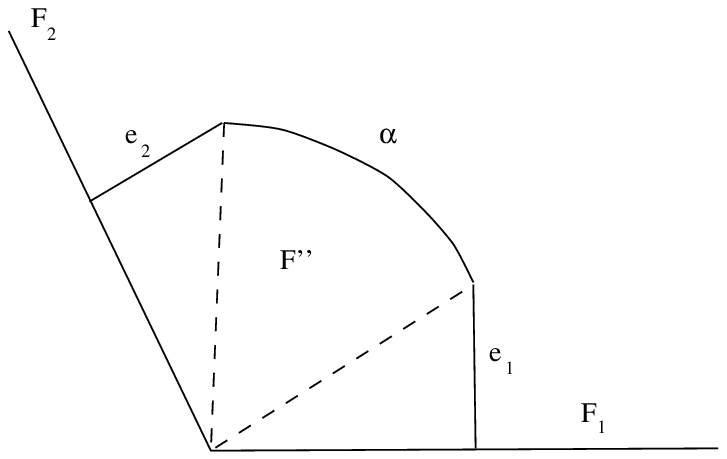|
\smallskip

More precisely, 
$F''$ is bounded by two triangular faces $F_1$ and
$F_2$ (in this collection of triangles),
two edges $e_1$ and $e_2$
of the Voronoi cell dual to the triangles,
and an arc $\alpha$
obtained from the truncation.  The two edges $e_1$ and
$e_2$ 
are perpendicular to the faces $F_1$ and $F_2$, 
respectively, by the
definition of Voronoi-Delaunay duality.  The edges $e_1$ and
$e_2$ meet
the faces $F_1$ and $F_2$, respectively,
by the construction of restricted cells.
The edges $e_1$ and $e_2$
cannot intersect at any point less than $\sqrt{2}$ from the
origin, because the point of intersection would be a point
equidistant from the four vertices of a simplex formed by
the vertices of $F_1$ and $F_2$.  
The simplex would have circumradius less than $\sqrt{2}$.
Its faces would then also have circumradius less
than $\sqrt{2}$, so
that the Delaunay simplex is small.  
This is impossible, since all small 
simplices have already been classified as
regions of the first type.

{\bf Lemma 4.4.}  In this context, assume that the faces $F_1$
and $F_2$ form an acute angle, and let $p$ be the point at which
the line through $e_1$ meets the plane through $F_2$.
Let $w_1$ and $w_2$ be the third vertices of 
the faces $F_1$ and $F_2$
respectively (that is, those other than $0$ and $v$).  If
the distance
from $0$ to $p$ is  at most $\sqrt{2}$,  then 
the simplex $(0,v,w_1,w_2)$ is small 
or a quasi-regular tetrahedron.

{\bf Proof:}  Suppose that the distance
from $p$ to $0$ is at most $\sqrt{2}$.  Let $c_1$
be the circumcenter of $F_1$. Let $S$ be the simplex
$(0,v,w_1,w_2)$.

We claim that $p$ lies in the interior of the triangle
$F_2$.  The bisecting line $\ell$ between $0$ and $v$
in the plane of $F_2$ contains $p$.  The line $\ell$
intersects two edges of $F_2$, once at $v/2$ and
once at some other point $p'$.  If $p$ is (strictly) outside
$F_2$, then $|p|>|p'|$.  This leads to a contradiction,
once we show $|p'|\ge\sqrt{2}$.  If $p'$ lies on
the edge between $v$ and $w_2$, this is clear, because
an elementary exercise shows that every point on the
line passing through $w_2$ and $v$ has distance at
least $\sqrt{2}$ from the origin.  Assume $p'$ lies
between $0$ and $w_2$, and consider $p'$ as a function
of $w_2$ and $v$.  Its length $|p'|$ attains its
minimum when $F_2$ is the right triangle $|v|=2$,
$|v-w_2|=2$, $|w_2|=2\sqrt{2}$.  Thus, $|p'|\ge\sqrt{2}$.

So $p$ lies in the interior of $F_2$.  The vertex 
$w_1$
has distance at most
$\sqrt{2}$ from $p$. No vertex $w_1$ can come
within $\sqrt{2}$ of an interior point of $F_2$ 
unless the circumradius of $F_2$
is at least $\sqrt{2}$.
If the circumradius of $F_2$ is $\sqrt{2}$, then $p$ is
the circumcenter of $S$, so that the circumradius of $S$
is $\sqrt{2}$, making $S$ a small simplex.
If the circumradius is greater than $\sqrt{2}$, then
since $|p-w_1|\le\sqrt{2}$, $p$ lies in the Voronoi cell at $w_1$.
Thus, $w_1$ clips (possibly degenerately) a standard region
across the faces $F_2$ from $w_1$, based at $0$, $v$, or $w_2$.
By Lemma 2.2 $(w_1,0,v,w_2)$ is a quasi-regular tetrahedron.\qed

We continue with our description of the figure in Diagram 4.3.
The arc $\alpha$ cannot be interrupted by
a further (Voronoi) edge of $F'$.  Such an edge would
be dual to a (Delaunay) face with vertices $0$, $v$, and some $v'$.
The circumradius of the triangle with these three vertices
would be less than $\sqrt{2}$ (because every edge in $F'$
comes within
distance $\sqrt{2}$ of the origin).  This contradicts the
construction of $F''$ with half-planes given above.
 This completes our discussion of the figure in Diagram 4.3.
We emphasize, however, that the edges $e_1$ or $e_2$
may degenerate to length $0$, and the circular arc $\alpha$ may
degenerate to a point.

This Voronoi face-wedge can be broken into three convex pieces:
the convex hull of $0$, $v/2$, and the circular arc,
and the convex hulls of $0$, $v/2$, 
and the edge $e_i$, for $i=1,2$.
The first piece has the third type, the others
have the fourth type.
The boundary condition
	$4/3\le b^2$ expresses the fact that the circumradius
	of a triangle with sides of length
at least $2$ cannot be less than $2\sqrt{3}/3$.
This completes the reduction to the four
given types.

Now we must show that each of the given types has density
at most $\delta_{oct}$.

{\bf Type 1:}  A small simplex that is not a
quasi-regular tetrahedron.
Let $S$ be a small simplex of Voronoi type
        with at least one edge longer than $2.51$.
        By the monotonicity properties of the circumradius,
        we know that the circumradius of $S$
         is at least $\rad(S(2,2,2,2.51,2,2))>1.3045$.
	Let $\delta(a,b,c)$ denote the density of the Rogers simplex
	$R(a,b,c)$ (see I.8.6).
        By Rogers's lemma (I.8.6.2), the six Rogers
        simplices $V_S$ have density less than $\doct$  and $\vor(S,V_S)<0$
        if the circumradius of the three faces is at least $1.207$
        ($\delta(1,1.207,1.3045)<\doct$). This condition on
        the circumradius of the faces holds whenever there
        are two edges longer than $2.51$ at the origin
        ($\eta(2.51,2,2)>1.207$) or whenever there are two oppositely
        arranged edges longer than $2.51$.  

Thus, to show that
        $\vor(S)<0$ for small simplices of Voronoi type,
        we must consider the following cases:  (1) one edge
        longer than $2.51$,
        (2) two adjacent edges longer than $2.51$, and (3)
        three edges longer than $2.51$ meeting at a vertex.
        These cases are covered by Calculation 4.5.2.
	In (1), we may assume that at least one of the conditions
	for compression scoring in Section 3.1 fails to hold.
        In Calculation 4.5.2.2, we may make the stronger assumption
        $\rad(S)<1.39$, for otherwise, the Rogers simplices
        at the origin have density at most
        $\delta(1,\eta(2,2,2.06),1.39)<\doct$ so that
        $\vor(S)<0$.

We rely on Calculation 4.5.1 for small simplices of compression type.
The appendix proves the result for simplices in an
        explicit neighborhood of $S(2\sqrt{2},2,2,2,2,2)$.
These calculations
        are established by methods of interval arithmetic described
        in Part I.  Source code appears in \cite{6}.

{\bf Calculation 4.5.1.}  If $S$ is a small
simplex that is not
a quasi-regular tetrahedron, then
$\Gamma(S)\le 0$.  If equality is attained, then
the simplex $S$ is congruent
to $S(2\sqrt{2},2,2,2,2,2)$ or to the simplex of zero
volume $S(2\sqrt{2},2,2,2\sqrt{2},2,2)$.

{\bf Calculation 4.5.2.} Assume $S$ is small.
 $\vor(S(y_1,\ldots,y_6))<0$ if $(y_1,\ldots,y_6)$ belongs
to any of the cells (1)--(11). Let $I$ denote the interval
$[2,2.51]$ and $L=[2.51,2\sqrt{2}]$.

{
\hbox{}
\obeylines
\parskip=0pt
\baselineskip=0.9\baselineskip

(1) \quad $L[2.06,2.51]I^4$,
(2) \quad $LI^2[2.06,2.51]I^2$,
(3) \quad $LI^3[2.08,2.51]I$,
(4) \quad $[2.06,2.51]I^2LI^2$,
(5) \quad $I[2.08,2.51]ILI^2$,
(6) \quad $I^3L[2.2,2.51]I$,
(7) \quad $I^3[2.58,2\sqrt{2}][2.12,2.51]I$,
(8) \quad $LI^3LI$,
(9) \quad $LI^3L^2$,
(10) \quad $I^3L^2I$,
(11) \quad $I^3L^3$.

}

{\bf Type 2:} The set $\{t x: 0\le t \le 1,\ x\in P_2\}$.
	In this case the density is increased by intersecting the
	set with a ball of radius $1.18$ centered
	at the origin.  The resulting
	intersection has density $1/1.18^2 <\delta_{oct}$, 
	as required.

{\bf Type 3:} The set $\{t x: 0\le t \le 1,\ x\in P_3\}$.
	The bounding circular arc of $P_3$ 
	has distance $\sqrt{2}$ from the origin.
	The set has the same density as a right circular cone,
	with base a disk of radius $\sqrt{2-h^2}$ and height $h$.
	This cone has volume $\pi(2-h^2)h/3$.  The solid angle at the
	apex of the cone is $2\pi(1-\cos\theta)$, where 
	$\cos\theta = h/\sqrt{2}$.
	This gives a density of $\sqrt{2}/(h^2+h\sqrt{2})$.
	This function is maximized
	over the interval $[1,1.18]$ at $h=1$.  The density
	is then at most $2-\sqrt{2} < \delta_{oct}$.

{\bf Type 4:} A Rogers simplex $R_1=R(a,b,\sqrt{2})$.
	where $1\le a\le 1.18$
	and $4/3\le b^2\le 2$.  

By Lemma I.8.6.2, the density of this simplex is at most
that of the
Rogers simplex $R_2=R(1,2\sqrt{3}/3,\sqrt{2})$.
This simplex has the density $\delta_{oct}$
of a regular octahedron.  (In fact, the regular octahedron may be
partitioned into simplices congruent to $R_2$ and its mirror.)
We see that the original simplex $R_1$
has density $\delta_{oct}$ exactly when,
in the notation of I.8.6.2,  $|{\bold s}_1|
= |{\bold s}_2|$, for all $\lambda_1$, $\lambda_2$, and $\lambda_3$
as above.  This implies that $a=1$ and $b=2\sqrt{3}/3$.  
This completes the proof of Theorem 4.1. \qed

{\bf Proposition 4.6.}  A cluster 
other than a quasi-regular tetrahedron attains
a score of $0\,\pt$ if and only if it is made up of
simplices
congruent to $S(2,2,2,2,2,2\sqrt{2})$, and possibly some additional
simplices of zero volume.

{\bf Proof.}
Types 2 and 3 always give strictly negative scores for regions
of positive volume.  According to Calculation 4.5, a region of
the first type with positive volume
gives a strictly negative score unless it is congruent
to $S(2,2,2,2,2,2\sqrt{2})$.  

Consider a region of the fourth type with score $0\,\pt$.
We must have $a=1$
and $b = 2\sqrt{3}/3$.
The circumradius of the faces $F_1$ and $F_2$ of Diagram 4.3
is then $2\sqrt{3}/3$.  This forces the faces $F_1$ and $F_2$
to be equilateral triangles of edge length $2$.  
The arc $\alpha$ in Diagram 4.3 must reduce to a point.
The edges $e_1$ and $e_2$ in Diagram 4.3 -- if they have positive
length -- must then meet at a point
at distance $\sqrt{2}$ from the origin.  This point
is a vertex of a Voronoi cell and the circumcenter of a Delaunay 
simplex $S$ (of circumradius $\sqrt{2}$).  The only simplex with two
equilateral faces of side 2 and $\rad(S) = \sqrt{2}$ is the
wedge of an octahedron $S=S(2,2,2,2,2,2\sqrt{2})$.
This is a small simplex.

The other possibility is that both the arc $\alpha$
and an edge (say $e_2$) degenerate
to length $0$.  In this case, Lemma 4.4 shows that the restricted
cell belongs to a small 
Delaunay simplex or a quasi-regular tetrahedron.  These
cases have already been treated.  \qed

\bigskip
\centerline{\bf References}
\bigskip

\baselineskip = 0.66\baselineskip
\parskip=0.5\baselineskip
\hbox{}

1.  Thomas C. Hales, 
	Remarks on the density of sphere packings in three dimensions,
	Combinatorica, 13 (2) (1993), 181--197.

2.  Thomas C. Hales, 
	The sphere packing problem, Journal of Comp. and Applied Math,
	44 (1992), 41--76.

3.  Thomas C. Hales,
	The status of the Kepler conjecture, Math. Intelligencer,
	(1994).

4.  Thomas C. Hales, 
	Sphere packings I, to appear in Discrete and Computational Geometry.

5.  Thomas C. Hales,
	Sphere packings III$_\alpha$, preprint.

6.  Thomas C. Hales,
	{\tt http://www.math.lsa.umich.edu/\~\relax
                hales/packings.html}

7.  W.-Y. Hsiang, On the sphere packing problem and the proof of 
	Kepler's conjecture, Int. J. of Math. 4, no.5, (1993), 
	739--831.

8.  C. A. Rogers, The packing of equal spheres, Proc. London Math.
	Soc. (3) 8 (1958), 609--620.

\vfill\eject
\centerline{\bf Appendix}
\bigskip

We give a direct argument that $\Gamma(S)\le 0\,\pt$, when the
lengths of 
a small simplex $S$ 
are within $0.001$ of $S_0=S(2\sqrt{2},2,2,2,2,2)$.
Set $S=S(y_1,y_2,y_3,y_4,y_5,y_6)$.
Write $y_1=2\sqrt{2}-f_1$, and $y_i = 2+f_i$, for $i>1$, where
$0\le f_i\le 0.001$.  Set $x_1 = y_1^2 = 8-e_1$ and
$x_i = y_i^2 = 4+e_i$, for $i>1$.  Then $0\le e_i\le 0.006$.
Recall from Section I.8.4 that
$$a(y_1,y_2,\ldots,y_6) = y_1y_2y_3 +
{1\over2}y_1(y_2^2+y_3^2-y_4^2) + {1\over2}y_2(y_1^2+y_3^2-y_5^2) +
{1\over2}y_3(y_1^2+y_2^2-y_6^2).
$$
Set
$$\align
a_0 &= a(y_1,y_2,y_3,y_4,y_5,y_6),\quad a_{00} = a(2\sqrt{2},2,2,2,2,2)=16+12\sqrt{2},\\
a_1 &= a(y_1,y_5,y_6,y_4,y_2,y_3),\quad a_{10} = a(2\sqrt{2},2,2,2,2,2)=16+12\sqrt{2},\\
a_2 &= a(y_4,y_2,y_6,y_1,y_5,y_3),\quad a_{20} = a(2,2,2,2\sqrt{2},2,2)=16,\\
a_3 &= a(y_4,y_5,y_3,y_1,y_2,y_6),\quad a_{30} = a(2,2,2,2\sqrt{2},2,2)=16.\\
\endalign
$$
Section I.8.4 
and the bounds on $f_i$ give $a_i\ge a_{i}^\minus$, where
$a_{0}^\minus=a_{1}^\minus=32.27$ and $a_{2}^\minus=a_{3}^\minus=15.3$.

Let $\Delta$ be the function of Section I.8.1, and
set $\Delta_0 = \Delta(8,4,4,4,4,4)$.
Set $t = \sqrt{\Delta(x_1,\ldots,x_6)}/2$, and $t_0 =
\sqrt{\Delta_0}/2$.  A simple calculus exercise shows
that 
$$\Delta(x_1,\ldots,x_6)\ge \Delta(8,4,4,x_4,4,4)=128-8e_4^2.$$
This gives $t\ge 5.628$.  Let $b_i = 2/(3(1+t_0^2/a_{i0}^2))$,
so that $b_0=b_1=(3+2\sqrt{2})/9$ and $b_2=b_3=16/27$.
Set $c_0 = -\delta_{oct}/6 + \sum_0^3 b_i/a_{i0} \approx -0.00679271$.
Then
  $$(t-t_0)c_0 = {(\Delta-\Delta_0)c_0\over 4(t+t_0)} \le
	{-2e_4^2 c_0\over (t+t_0)} \le 0.002 e_4^2 < 0.0002 f_4.$$
We are ready to estimate $\Gamma(S)$.  An argument parallel to
that of Lemma I.9.1.1 gives
$$\Gamma(S) \le \Gamma(S_0) + (t-t_0) c_0 + 
	t \sum_{i=0}^3 {b_i (a_{i0}-a_i)\over a_{i0}^2} +
	t \sum_{i=0}^3 {b_i (a_{i0}-a_i)^2\over a_{i0}^2 a_{i}^\minus}.\tag1
$$
The two sums on the right-hand side are polynomials in $f_i$ with
no constant terms.  To give an upper bound on these polynomials,
write them as a sum of monomials,
and discard the negative monomials of order greater than 2.
The positive monomials of order greater than 2 are dominated
by $$f_1^{d_1} f_2^{d_2}\cdots f_6^{d_6} 
	\le (0.001)^{d_1+\cdots +d_6-1}
	(f_1+f_2+\cdots +f_6).$$
This approximation shows that the first sum in (1)
is at most $-0.005 f_1 - 0.04f_4-0.03 (f_2+f_3+f_5+f_6)$ and the
second sum in (1) is at most
$0.00056 (f_1+f_2+\cdots+f_6)$.
The result easily follows.

This argument is easily adapted to a neighborhood of $S_1 =
S(2\sqrt{2},2,2,2\sqrt{2},2,2)$.  In this case, for
$i=1,\ldots,4$, we have $a_{i0}= 16+8\sqrt{2}$, $b_i = 2/3$,
$a_i^\minus=27$,
$t_0=0$, $c_0\approx -0.0225$, and $t\ge0$.  A similar
argument leads to
the conclusion that $\Gamma(S)<\Gamma(S_i)=0\,\pt$, if $S$
is a small simplex such that
$S\ne S_1$, and the lengths of the edges of $S$ are within
$0.01$ of those of $S_1$.

\bigskip
\hbox{}
\smallskip

Research supported by the NSF.

Dept. of Math, University of Michigan

Version: 4/12/95, 3/28/96.

\bye